\documentclass[twoside,notitlepage,11pt]{article}

\usepackage{amssymb}
\usepackage[leqno]{amsmath}
\usepackage{amsfonts}
\usepackage{amsopn}
\usepackage{amstext}
\usepackage{amsthm}

\providecommand{\cal}{\mathcal}
\renewcommand{\Bbb}{\mathbb}

\newenvironment{pf}{\begin{proof}}{\end{proof}}

\newcommand{\Ha}{{\cal{H}}}

\newcommand{\Zee}{{\Bbb{Z}}}

\newcommand{\lam}{{\lambda}}
\newcommand{\al}{\alpha}

\renewcommand{\phi}{\varphi}
\renewcommand{\rho}{\varrho}

\newcommand{\rest}{\restriction}

\newcommand{\loe}{\leqslant}
\newcommand{\goe}{\geqslant}

\newcommand{\subs}{\subseteq}

\newcommand{\id}{\operatorname{id}}

\newcommand{\liminv}{\varprojlim}
\newcommand{\by}{/}

\newcommand{\setof}[2]{\{#1\colon #2\}}

\newcommand{\seq}[1]{\langle #1 \rangle}

\newcommand{\sett}[2]{\{#1\}_{#2}}
\newcommand{\map}[3]{#1\colon #2 \to #3} % A function
\newcommand{\suppt}{\operatorname{suppt}}

\newcommand{\Sig}{\Sigma}
\newcommand{\R}{\ensuremath{\mathcal R}}

\newcommand{\invsys}[5]{\langle {#1}_{#4},{#2}_{#4}^{#5},#3 \rangle}

\renewcommand{\S}{\mathbb S}
\newcommand{\T}{\mathbb T}
\newcommand{\rk}{\operatorname{rk}}

\newcommand{\cube}[1]{{[0,1]^{#1}}}

\newcommand{\Hone}{\ensuremath{H^1}} % the first cohomology functor
\newtheorem{tw}{Theorem}[section]
\newtheorem{wn}[tw]{Corollary}
\newtheorem{lm}[tw]{Lemma}

\theoremstyle{definition}

\theoremstyle{remark}

\newtheorem{question}{Question}

\title{Valdivia compact Abelian groups}
\author{
{\sc Wies{\l}aw Kubi\'s}\\
{\small Instytut Matematyki, Akademia \'Swi\c{e}tokrzyska, Kielce, Poland}
\\
\texttt{wkubis@pu.kielce.pl}
}

\begin{document}
\maketitle

\begin{abstract} 
Let \R\ denote the smallest class of compact spaces containing all metric compacta and closed under limits of continuous inverse sequences of retractions. Class \R\ is striclty larger than the class of Valdivia compact spaces.
We show that every compact connected Abelian group which is a topological retract of a space from class \R\ is necessarily isomorphic to a product of metric compacta. This completes the result of \cite{KU}, where a compact Abelian group outside class \R\ has been described.

\ 

\noindent {\bf MSC (2000):}
Primary: 54D30. Secondary: 54C15, 22C05.

\noindent {\bf Keywords:} Compact Abelian group, Valdivia compact, class \R, product, retraction, first cohomology functor.
\end{abstract}

\section{Introduction}

In the last years there has been some significant interest in the theory of {\em Valdivia compact spaces} \cite{AMN, DG}, i.e. those compact spaces $K$ which are embeddable into a Tikhonov cube $\cube\kappa$ in such a way that $K\cap \Sigma(\kappa)$ is dense in $K$, where $\Sigma(\kappa)$ denotes the $\Sigma$-product of $\kappa$ copies of $[0,1]$. The main motivation comes from the fact that every Valdivia compact $K$ has a resolution consisting of retractions onto Valdivia compacta of smaller weights. Consequently, the Banach space $C(K)$ has a {\em projectional resolution of the identity}, a useful property with several important consequences. For details we refer to \cite{BKT}, \cite{DG}, \cite[Chapter VI]{DGZ} and \cite[Chapter 6]{Fa}.
Several interesting results on Valdivia compacta were proved by Ond\v rej Kalenda \cite{Kalenda_segment, Kalenda4, Kalenda5, Ka_images} (see also his survey article \cite{Ka_survey}). 

One of the questions, left open in \cite{Ka_survey}, asked whether the class of Valdivia compacta was stable under open images. Typical examples of open surjections are epimorphisms of compact groups. Since every compact group is an epimorphic image of a product of metric compact groups and every product of metric compacta is Valdivia compact, it is natural to ask whether all compact groups are Valdivia compact.
A counterexample has been found by Vladimir Uspenskij and the author \cite{KU} -- a compact connected Abelian group whose Pontryagin dual is indecomposable, i.e. not representable as a direct sum of two proper subgroups.

The purpose of this note is to complete the result of \cite{KU}. Namely, we show that every connected Abelian group $G$ which is Valdivia compact must be isomorphic (as a topological group) to a product of metrizable compact groups. In fact, we prove a stronger result involving class \R, introduced in \cite{BKT} and defined to be the smallest class of spaces containing all metric compacta and closed under limits of continuous inverse sequences whose bonding maps are retractions. The study of class \R\ has the same functional-analytic motivations as for the (strictly smaller) class of Valdivia compacta, see \cite{BKT} and \cite{K_classR}. It has been proved in \cite{KU} that the mentioned above Abelian compact group does not belong to class \R.

Using similar methods as in \cite{KU}, we show that a compact connected Abelian group which is a topological retract of some space from class \R\ must be isomorphic to a product of metric groups and consequently is Valdivia compact.
We also deduce a similar statement for disconnected Abelian groups, using a well known topological product decomposition.

Let us admit that it is not known whether class \R\ is stable under retractions.

\section{Preliminaries}

A {\em topological group} is a group together with a Hausdorff topology for which the group operations are continuous. Given a locally compact Abelian group $G$, the {\em Pontryagin dual} $\hat G$ of $G$ is the group of all continuous homomorphisms $\map \chi G\T$ endowed with the pointwise convergence topology, where $\T$ is the circle group. It is well known that the Pontryagin dual of the Pontryagin dual of an Abelian group $G$ is isomorphic to $G$ and $\hat G$ is compact (and connected) iff $G$ is discrete (and torsion-free). For more information on Pontryagin duality we refer to \cite{HofMor} and \cite{HewRos}.

All topological spaces under consideration are assumed to be compact Hausdorff. 
By a {\em retraction} we mean a continuous map $\map fXY$ which is {\em right-invertible}, i.e. $fj=\id_Y$ for some continuous map $\map jYX$.
A {\em retractive inverse sequence} is a continuous inverse sequence $\S=\invsys Kr\delta\al\beta$ in which all bonding maps $r^\beta_\al$ are retractions (equivalently: each $r^{\al+1}_\al$ is a retraction). Recall that a sequence $\S$ is {\em continuous} if $K_\rho$ together with maps $\sett{r^\rho_\al}{\al<\rho}$ is the limit of $\S\rest\rho = \invsys Kr\rho\xi\eta$ for every limit ordinal $\rho<\delta$.

Recall that a compact space $K$ is {\em Valdivia compact} if $K$ can be embedded into a Tikhonov cube $\cube \kappa$ so that $K\cap \Sig(\kappa)$ is dense in $K$, where $\Sigma(\kappa)=\setof{x\in \cube\kappa}{|\suppt(x)|\loe\aleph_0}$ and $\suppt(x)=\setof{\al}{x(\al)\ne0}$.
Class \R\ is defined to be the smallest class of spaces containing all metric compacta and closed under limits of retractive inverse sequences. It is well known (see e.g. \cite{K_classR}) that class \R\ is strictly larger than the class of Valdivia compacta.
There is a natural ordinal rank $\rk_\R$ on class \R. Namely, we define $\rk_\R(K)=0$ iff $K$ is a metric compact and we declare $\rk_\R(K)\loe \beta$ if $K=\liminv\S$, where $\S=\invsys Kr\kappa\xi\eta$ is a retractive sequence such that $\rk_\R(K_\xi)<\beta$ for every $\xi<\kappa$. Finally, we set $\rk_\R(K)=\beta$ if $\rk_\R(K)\loe\beta$ and $\rk_\R(K)\not\loe \al$ for any $\al<\beta$.

\section{Continuous functors and Abelian groups}

The key tool needed to obtain our main result is a continuous contravariant functor from the category of compact spaces to the category of discrete Abelian groups. Namely, as in \cite{KU}, we use the first cohomology functor $\Hone$, in the sense of the theory of sheaves \cite{Go} or the \v Cech theory \cite[Chapter 6]{Sp}. More precisely, we let $\Hone(X)=\Hone(X,\Zee)$, i.e. we use the integers as the coefficient group. Functor $\Hone$
has the following properties:
\begin{enumerate}
	\item[(1)] $\Hone(X)$ is countable whenever $X$ is second countable. 
	\item[(2)] $\Hone(f)$ is a monomorphism whenever $f$ is a continuous surjection.
	\item[(3)] If $X=\liminv\invsys Xp\kappa\xi\eta$ then $\Hone(X)$ is the inductive limit of the sequence of groups $\seq{\Hone(X_\xi),\Hone(p^\eta_\xi),\kappa}$.
\end{enumerate}
For the proofs we refer to \cite{Sp} or \cite{KU}. Another important feature of the above functor is that $\Hone(G)$ coincides with the Pontryagin dual of $G$, whenever $G$ is a compact connected Abelian group (see \cite{KU} for more details).

\begin{tw}\label{woegfijaif} Assume that $X$ is a retract of a space from class $\R$. Then $\Hone(X)$ is isomorphic to a direct sum of countable groups. \end{tw}

For the proof we need two statements about discrete Abelian groups\footnote{Although these statements seem to be well known, we could not find any reference, therefore we supply the proofs.}. Let us denote by $\Ha$ the class of all Abelian groups which are isomorphic to a direct sum of countable groups. It turns out that $\Ha$ has a very simple structure.

\begin{lm}\label{ijegwpjgp} Assume that $\sett{G_\al}{\al<\kappa}$ is a continuous increasing chain of subgroups of an Abelian group $G$ such that $G=\bigcup_{\al<\kappa}G_\al$, $G_0=0$ and $G_{\al+1}=G_\al\oplus H_\al$ for every $\al<\kappa$. Then
$$G=\bigoplus_{\al<\kappa}H_\al.$$
In particular, if $\setof{H_\al}{\al<\kappa}\subs\Ha$ then also $G\in\Ha$.
\end{lm}

\begin{pf} Denote by $G'$ the algebraic sum $\sum_{\al<\kappa}H_\al\subs G$, i.e. $G'$ consists of all elements of the form $x_0+\dots+x_{k-1}$, where $x_i\in H_{\al_i}$, $i<k$. Then $G_0\subs G'$ and $G_\xi\subs G'$ implies $G_{\xi+1}\subs G'$. Thus, by induction and by the continuity of the chain, we deduce that $G'=G$.

Now suppose $x_0+\dots+x_k=0$, where $x_i\in H_{\al_i}$ for $i\loe k$. Assume  $0<\al_0<\al_1<\dots<\al_k$ and $\al_i$ is a minimal ordinal $\al$ such that $x_i\in H_\al$ ($i\loe k$). By the continuity of the chain, we have $\al_k=\xi+1$ and hence $x_0+\dots+x_{k-1}\in G_\xi$ and $x_k\in H_\xi$. Thus $x_k=0$, because $G_\xi\cap H_\xi=0$. This shows that $H_\eta\cap \sum_{\al\ne\eta}H_\al=0$ for every $\eta<\kappa$.
\end{pf}

\begin{lm}\label{iojwetttepw} Class $\Ha$ is closed under direct summands. That is, if $G=H\oplus K\in\Ha$ then $H,K\in\Ha$.
\end{lm}

\begin{pf} We use induction on the cardinality of the group. Of course, the claim is true for countable groups. Fix $\kappa>\aleph_0$ and suppose that the statement holds for groups of cardinality $<\kappa$. Fix $G=\bigoplus_{\al<\kappa}G_\al$, where each $G_\al$ is a countable Abelian group. Let $H$ be a direct summand of $G$ and let $\map hGH$ be a group epimorphism such that $h\rest H=\id_H$. Fix a cardinal $\chi$ big enough so that $h\in H(\chi)$ and $G\subs H(\chi)$. Fix an increasing continuous chain $\sett{M_\al}{\al<\kappa}$ of elementary substructures of $H(\chi)$ such that $h\in M_0$, $H\subs\bigcup_{\al<\kappa}M_\al$ and $|M_\al|<\kappa$ for every $\al<\kappa$.
Let $H_\al=H\cap M_\al$. Then $\sett{H_\al}{\al<\kappa}$ is an increasing continuous chain of subgroups of $H$ and $H=\bigcup_{\al<\kappa}H_\al$. 

Fix $\al<\kappa$ and let $M=M_\al$. Observe that $$G\cap M=\bigoplus_{\xi\in \kappa\cap M}G_\al.$$
Indeed, if $\xi\in\kappa\cap M$ then $G_\xi\subs M$, because $|G_\xi|\loe\aleph_0$. Thus $\bigoplus_{\xi\in \kappa\cap M}G_\xi\subs M$. On the other hand, if $x\in G\cap M$ and $x=x_0+\dots+x_{k-1}$, where $x_i\in G_{\xi_i}$, then $\setof{\xi_i}{i<k}\subs M$ and hence $x\in \bigoplus_{\xi\in \kappa\cap M}G_\al$.
Let $\map{p_\al}G{G\cap M}$ denote the canonical projection. By elementarity, $g_\al=hp_\al$ is a homomorphism of $G$ onto $H_\al$ which is identity on $H_\al\subs G\cap M$. Thus $H_{\al+1}=H_\al\oplus K_\al$, where $K_\al=H_{\al+1}\cap \ker(g_\al)$. Applying Lemma \ref{ijegwpjgp}, we get 
$$H=H_0\oplus\bigoplus_{\al<\kappa}K_\al.$$
By inductive hypothesis, $H_0,K_\al\in\Ha$, because both of these groups are direct summands of $G\cap M_{\al+1}$ and $|G\cap M_{\al+1}|<\kappa$.
Hence $H\in\Ha$, which completes the proof.
\end{pf}

\begin{pf}[Proof of Theorem \ref{woegfijaif}] By Lemma \ref{iojwetttepw}, it suffices to show that $\Hone(X)\in\Ha$ whenever $X\in\R$. We use induction on $\rk_\R$. By property (1) of the functor $\Hone$, the statement is true for spaces of \R-rank $0$. Fix an ordinal $\beta>0$ and assume $\Hone(X)\in\Ha$ whenever $\rk_\R(X)<\beta$. Fix $X\in \R$ with $\rk_\R(X)=\beta$ and let $\S=\invsys Xr\kappa\xi\eta$ be a continuous retractive sequence with $X=\liminv\S$ and such that $\rk_\R(X_\xi)<\beta$ for every $\xi<\kappa$. Let $G=\Hone(X)$ and let $G_\xi=\Hone(X_\xi)$. By properties (2) and (3), we may identify each $G_\xi$ with a subgroup of $G$ so that $\sett{G_\xi}{\xi<\kappa}$ becomes a continuous increasing chain with $G=\bigcup_{\xi<\kappa}G_\xi$. Further, each $G_\xi$ is a direct summand of $G$. By inductive hypothesis $\setof{G_\xi}{\xi<\kappa}\subs \Ha$. By Lemma \ref{ijegwpjgp}, we have that $G=G_0\oplus\bigoplus_{\xi<\kappa}H_\xi$, where $H_\xi$ is such that $G_{\xi+1}=G_\xi\oplus H_\xi$ for every $\xi<\kappa$. Since $H_\xi$ is a direct summand of $G_\xi$, we deduce, using Lemma \ref{iojwetttepw}, that $\setof{H_\xi}{\xi<\kappa}\subs\Ha$. Thus $G\in\Ha$. \end{pf}

\begin{tw} Assume that $G$ is a compact connected Abelian group which is at the same time a topological retract of some space from class \R. Then $G$ is isomorphic, in the category of topological groups, to a product of metrizable compact groups.
\end{tw}

\begin{pf}
Let $H$ denote the Pontryagin dual of $G$. Since $G$ is connected, $H$ is discrete and torsion-free, therefore isomorphic to $\Hone(G)$ (see \cite[Proposition 2.5]{KU} for a proof of this well known fact).
By Theorem \ref{woegfijaif}, $H$ can be decomposed into a direct sum of countable groups. Thus $G$ is a product of metric groups, because $G$ is the dual of $H$ and Pontryagin duality turns direct sums into products. \end{pf}

It is well known that every compact group is a Dugundji space (i.e. an absolute extensor for the class of $0$-dimensional compacta). It has been proved in \cite{KM} that $0$-dimensional Dugundji compacta (equivalently: retracts of Cantor cubes) are Valdivia compact. In the case of compact groups this result becomes trivial: 
every infinite $0$-dimensional compact group is homeomorphic to a Cantor cube $2^\kappa$ for some cardinal $\kappa\goe\aleph_0$, see \cite[Theorem 9.15]{HewRos}. In general, every compact group $G$ is homeomorphic to $G_0\times H$, where $G_0$ is the component of the identity and $H$ is a $0$-dimensional group isomorphic to $G\by G_0$, see \cite[Corollary 10.37]{HofMor}.
%In particular, $H$ is either finite or homeomorphic to a Cantor cube. 
Thus, using the above theorem, we conclude the following:

\begin{wn} Let $G$ be a compact Abelian group. The following properties are equivalent:
\begin{enumerate}
	\item[(a)] $G$ is a topological retract of some space from class \R.
	\item[(b)] $G$ is Valdivia compact.
	\item[(c)] $G$ is homeomorphic to a product of metric compacta.
	\item[(d)] The identity component of $G$ is isomorphic to $\prod_{\xi<\lam}H_\xi$, where $\lam$ is a cardinal and $H_\xi$ is a compact metric group for every $\xi<\lam$.
\end{enumerate}
\end{wn}

We finish with the following natural

\begin{question}
Does there exist a (non-commutative) Valdivia compact group $G$ which, as a topological space, is not homeomorphic to any product of metric compacta?
\end{question}

\end{document}